\documentclass[freqn,leqno,10pt]{amsart}
\usepackage{amsthm} 
\usepackage{graphicx}
\usepackage{array}
\usepackage{amsmath}
\usepackage{amssymb}
\usepackage{amsfonts}
\usepackage{latexsym}
\newtheorem{theorem}{Theorem}
\newtheorem{definition}{Definition}

\newtheorem{proposition}{Proposition}

\newtheorem{corollary}{Corollary}
\newtheorem{remark}{Remark}
\newtheorem{fact}{Fact}

\newtheorem{example}{Example}

\begin{document}
\title{Seifert matrices of Brunnian links}
\author{Maki Nagura}
\maketitle
\begin{center}{\footnotesize{\it Department of Applied Mathematics, Yokohama National University, \\
79-5 Tokiwadai, Hodogaya, Yokohama 240-8501, Japan} \\
maki@ynu.ac.jp}
\end{center}

\begin{abstract} 
We show some properties of a Seifert matrix of an $n$-component Brunnian link. In particular, we give a necessary and sufficient condition for a matrix to be a Seifert matrix of a $2$-component Brunnian link up to S-equivalence. 
\end{abstract}

\bigskip

\section{Introduction} 
\medskip

Throughout this paper, we shall work in the PL category and study tame oriented links in the 3-sphere ${\mathbb S}^3$ that are called Brunnian links. Brunnian links are constructed by Hermann Brunn in his paper ~\cite{1} and are defined as follows:

\begin{definition}(cf. ~\cite{7}) If a link is non-trivial and every proper sublink is trivial, then we say that it has the Brunnian property, or it is Brunnian.  
\end{definition}

For example, the Borromean rings, the Whitehead's link, and the Milnor links are Brunnian. These links are algebraically split, but in this paper, we shall treat all Brunnian links.
In section 2, we introduce a C-complex spanning a Brunnian link, and in section 3, we describe a method for obtaining a Seifert surface of a Brunnian link. In section 4, we show properties of a Seifert matrix of a Brunnian link. In section 5, we prove a necessary (Theorem \ref{necess}) and sufficient (Theorem \ref{suff}) condition for a matrix to be a Seifert matrix of a $2$-component Brunnian link up to S-equivalence(cf. \cite{8}), and show some examples of alternations. Finally, we show a Seifert matrix of an $n$-component Brunnian link, $n\geqq 3$, in section 6.

\bigskip

\section{Preliminaries}
\medskip

To construct a Seifert surface for a Brunnian link, we need the following properties: 

%\begin{remark} Let $L$ be an n-component Brunnian link. Then, there is a diagra%m for $L$ such that $(n-1)$ components of $L$ are all disjoint.
%\end{remark}

\begin{fact} Let $L=\cup_{i=1}^n K_i$ be an n-component Brunnian link. Then, for any $i=1, 2, \ldots, n$, there is a diagram for $L$ such that components $K_j$ are all disjoint, where $j\in\{1, 2, \ldots , n\}\setminus\{i\}$.
\end{fact}

\begin{proposition} \label{pro1}Let $L=\cup_{i=1}^n K_i$ be an $n$-component Brunnian link. Then, there are (2-)disks $D^2_1, D^2_2, \ldots , D^2_n$ satisfying the following $(i)$ and $(ii)$:
\begin{enumerate}
\item [$(i)$] $\partial D^2_i=K_i$, $i=1, 2, \ldots, n$, and ;
\item [$(ii)$] $D^2_i\cap D^2_j=\emptyset$, $i\neq j, \,i, \,j\in\{1, 2, \ldots, n-1\}$. 
\end{enumerate}
\end{proposition}
\begin{proof} Choose a diagram for $L$ so that components $K_1, K_2\ldots, K_{n-1}$ are all disjoint circles by Fact 1. And then, span a disk $D^2_i$ to each $K_i$, $i=1, 2, \ldots, n-1$ so that $D^2_i\cap D^2_j=\emptyset$ for any $j\in\{1, 2, \ldots, n-1\}\setminus \{i\}$. Since $K_n$ is unknotted, it spans (in the 3-sphere ${\mathbb S}^3$) a disk without self-intersections. Then, the disks $D^2_i$, $i=1, 2, \ldots, n$, satisfy $(i)$ and $(ii)$.
\end{proof}

\begin{corollary} Let $L=\cup_{i=1}^n K_i$ be an $n$-component Brunnian link. Then, there is a Seifert complex $\cup _{i=1}^n D^2_i$ satisfying $(i)$ and $(ii)$ in Proposition \ref{pro1}, where a Seifert complex is a finite union of compact pl embedded surfaces that are in general position in ${\mathbb S}^3$ and intersect transversely, and $D^2_i$ is a disk.
\end{corollary}

\begin{proof} By Proposition 1, there are disks $D^2_1, D^2_2, \ldots , D^2_n$ satisfying $(i)$ and $(ii)$ in Proposition 1. Make each pair $D^2_i$ and $D^2_n$, $i=1, 2, \ldots , n-1$, of the disks transverse and arrange $D^2_1, D^2_2, \ldots , D^2_n$ in general position in ${\mathbb S}^3$, then the union $\cup_{i=1}^n D^2_i$ is a Seifert complex satisfying $(i)$ and $(ii)$ in Proposition \ref{pro1}. 
\end{proof}

Let $\cup_{i\in I}S_i$ be a Seifert complex, where $I$ is a finite set. Then a (connected) component of an intersection $\cap_{j\in J} S_j$, if it exists, is called a singularity of the Seifert complex $\cup_{i\in I}S_i$, where $J\subset I$. A component of a Seifert complex is a subset of the Seifert complex whose preimage is a (connected) component. Note that a Seifert complex $S=S_1\cup S_2$ that consists of two components has singularities of tree types: clasp (or C), ribbon (or R), and circle, where $S_1$ and $S_2$ are compact pl embedded surfaces in ${\mathbb S}^3$. The details are referred to \cite{15}, \cite{16}.  
\smallskip

A Seifert complex is called a C-complex if all singularities are clasp, and is called an RC-complex if all singularities are ribbon or clasp (cf. \cite{3}, \cite{4}).
\smallskip

The following result is proved by C.Cooper.

\begin{proposition} \cite{4}
Any pair of Seifert surfaces for a link may be isotoped keeping their boundaries fixed to give a C-complex.
\end{proposition}

%Therefore, the following proposition is trivial.

\begin{proposition} Let $L=\cup_{i=1}^n K_i$ be an $n$-component Brunnian link. Then, there is an RC-complex (and a C-complex) $\cup _{i=1}^n D^2_i$ satisfying $(i)$ and $(ii)$ in Proposition \ref{pro1}, where each $D^2_i$ is a disk.
\end{proposition}

The following proof is referred to \cite{4}.

\begin{proof} By Corollary 1, there is a Seifert complex $\cup _{i=1}^n D^2_i$ satisfying $(i)$ and $(ii)$ in Proposition \ref{pro1}. First, remove an outermost circle singularity, if it exists, on $D^2_i\cap D^2_n\subset D^2_i$ for some $i\in\{1, 2, \ldots, n-1\}$, by pushing $D^2_n$ along an arc going from $\partial D^2_n$ to that circle singularity. Note that $D^2_n$ has no self-intersections by its construction. This transforms the circle singularity into a ribbon singularity. Continue in this way, until all circle singularities on $\cup_{i=1}^n D^2_i$ are removed. Then, the resulting Seifert complex, bounds $L=\cup_{i=1}^n K_i$, is an RC-complex. 

Next, remove a ribbon singularity, if it exists, on $D^2_j\cap D^2_n\subset D^2_j$ for some $j\in\{1, 2, \ldots, n-1\}$, by pushing $D^2_j$ along an arc going from $\partial D^2_j$ to the ribbon singularity to replace it by two clasps. Continue in this way, until all ribbon singularities on $\cup_{i=1}^n D^2_i$ are removed. Then, the resulting Seifert complex, still bounded $L=\cup_{i=1}^n K_i$, is a C-complex.
\end{proof}

\bigskip

\section{Constructing a Seifert surface for a Brunnian link}
\medskip

\noindent {\bf Algorithm.} Let $L$ be an $n$-component Brunnian link. First, span $L$ by a union of disks satisfying $(i)$ and $(ii)$ in Proposition \ref{pro1}, and then transform the union into a Seifert complex as in Corollary 1, and then transform the Seifert complex into a C-complex, say $C_L$ as in Proposition 3. Finally, remove each singularity on $C_L$, if it exists, by performing an {\it orientation preserving cut} along the singularity (see ~\cite{5}). So we obtain a Seifert surface for $L$. 
\bigskip

To obtain a Seifert matrix for an $n$-component Brunnian link $L$, we shall use a Seifert surface constructed by the above \lq\lq algorithm\rq\rq\, starting from a C-complex $C_L=\cup_{i=1}^n D^2_i$ spanning $L$, as in Proposition 3. We note that double points on $C_L$ are clasp. Let $C_L$ and $S_L$ be a Seifert complex and a Seifert surface for $L$ obtained by the above algorithm. Then the surface in $S_L$ that is deformed from $D_i^2$ is denoted by $V_i$ (see Figure ) and we denote a disk as in Figure by $D_i$. 
\bigskip

\section{Theorems and Properties}
\medskip

In this section, we denote the following $n\times (n-1)$ matrix by $F_{n}=(f_{i,j})$, $i\in\{1 ,\ldots ,n\}$, $j\in\{1 ,\ldots ,n-1\}$, $n\in{\mathbb N}$: if $i-j=0$ or $1$, then $f_{i,j}=1$, otherwise, $f_{i,j}=0$. A diagonal matrix whose diagonal element is $1$ or $-1$ is denoted by $E$.

In the next theorem, we show how to take an oriented base of the first homology group on a Seifert surface for a 2-component Brunnian link obtained in section 3, and construct a Seifert matrix by using the oriented base.

\begin{theorem}\label{Smatrix} For any 2-component Brunnian link $L$, there is a Seifert matrix $M_L$ for $L$ of the following form:
$$
M_L=\left(
\begin{array}{c|c}
E & F_{n}  \\ \hline
0 & H \\
\end{array}
\right).
$$
\end{theorem}

\begin{proof} Let $L=K_1\cup K_2$. Let $C_L=D^2_1\cup D^2_2$ be a C-complex and $S_L$ be a Seifert surface for $L$ obtained according to \lq\lq Algorithm\rq\rq\, in section 3. As a set of oriented loops representing a base of $H_1(\textrm{Int}(S_L))$, we may take ${\mathcal L}\equiv\{a_1, a_2, \ldots, a_{n}, b_{1}, b_2, \ldots, b_{n-1}\}$ such that
$$
\left\{
\begin{array}{@{\,}ll}
D_1\cap (a_k\cap b_j)=\emptyset\, , \\

D_1\cap (b_i\cap b_j)=\emptyset\, ,  \\

D_2\cap (a_k\cap b_j)=\textrm{at most one point}\, , \\

D_2\cap (b_i\cap b_j)=\textrm{at most one point}\,,    \\
\end{array}
\right.
$$
$$
a_k\cap b_j=\left\{
\begin{array}{@{\,}ll}
\textrm{ one point } & \textrm{if} \,\, \, k=j, \, j-1,  \\

\, \, \emptyset & \textrm{if}  \,\, \, k\neq j, \, j-1, \,\textrm{and}\,;
\end{array}
\right.
$$
\noindent $lk(a_i, b_i^+)=lk(a_{i+1}, b_i^+)=1$ for any $i$, where an integer $n$ is the number of clasp singularities on the disk $D^2_1$ of $C_L$ and $k\in\{1, 2, \ldots , n\}$, $i,j\in\{1, 2, \ldots , n-1\}$ (see Figures \ref{cbase}-\ref{aibi}). 

\begin{figure}
\begin{center}
\includegraphics[clip]{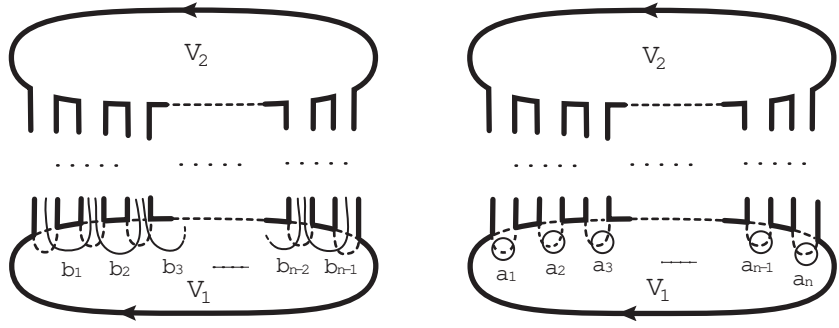}
\end{center}
\caption{choice of bases}
\label{cbase}
\end{figure}

\begin{figure}
\begin{center}
\includegraphics[clip]{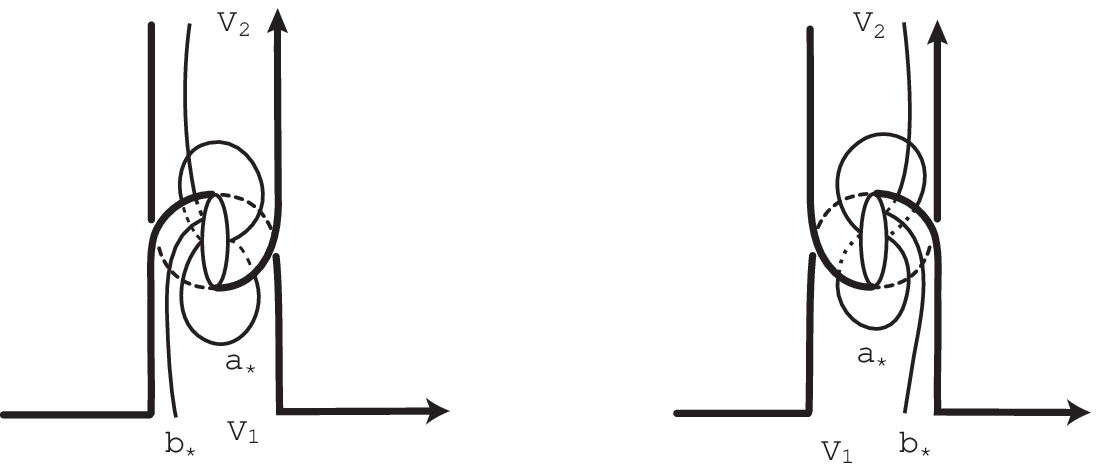}
\end{center}
\caption{}
\label{aibi}
\end{figure}

By using the base ${\mathcal L}$ of $H_1(\textrm{Int}(S_L))$, we obtain a Seifert matrix $M=(m_{ij})$, $i,j\in\{1, 2, \ldots, 2n-1\}$ for $L$ as follows:  
$$
m_{ij}=\left\{
\begin{array}{@{\,}ll}
\smallskip

lk(a_i, a_j^{+}) & \textrm{if} \,\,\, i, j\in\{ 1, 2, \ldots , n\}, \\ 
\smallskip

lk(a_i, b_{j-n}^{+}) & \textrm{if} \,\,\, i\in\{ 1, 2, \ldots , n\},\, j\in\{ n+1, n+2, \ldots, 2n-1\}, \\
\smallskip

lk(b_i, a_j^{+}) & \textrm{if} \,\,\, i\in\{ n+1, n+2, \ldots , 2n-1\},\, j\in\{ 1, 2, \ldots , n\},\\
\smallskip

lk(b_{i-n}, b_{j-n}^{+}) & \textrm{if} \,\,\, i, j\in\{ n+1, n+2, \ldots , 2n-1\},\\
\end{array}
\right.
$$
\noindent where $x^{+}$ means a loop obtained by lifting the loop $x$ over $S_L^{+}$ slightly and $X^{+}$ means the positive side of the surface $X$. Then, the Seifert matrix $M$'s form is $M_L$ as in Theorem 1.

\end{proof}

\begin{remark}
Given a Seifert matrix $M_L$ of a 2-component Brunnian link $L$ as in the proof of Theorem \ref{Smatrix}, a matrix which is S-equivalent to $M_L$ might not be a Seifert matrix of $L$.
\end{remark}
\medskip

From now on, we take a Seifert surface for a 2-component Brunnian link and an oriented base $\{a_1, a_2, \ldots, a_{n}, b_{1}, b_2, \ldots, b_{n-1}\}$ of the first homology group on the surface as in Figure \ref{cbase}.
\medskip

\begin{remark} If a 2-component Brunnian link $L$ is algebraically split, then the submatrix $E=(e_{i,j})$ of a Seifert matrix $M_L$ as in Theorem \ref{Smatrix} satisfies the equality: $\sum_{i=1}^n e_{ii}=0$.
\end{remark}
\medskip

To describe properties of a submatrix $H=(h_{i,j})$ obtained by our choice of bases (see Figure \ref{cbase}), we need the following notations.

\begin{definition} Let $\{b_1, b_2 ,\ldots , b_{n-1}\}$ be a set of oriented loops as in Figure \ref{cbase} and let $H=(h_{i,j})=\left(lk(b_{i}, b_{j}^{+})\right)$, $i, j\in\{ 1, 2, \ldots , n-1\}$. Then we call $h_{i,j}-h_{j,i}$ the alternation of $h_{i,j}$, denoted by $h_{i,j}^A$, where $i<j$, and we define if $h_{i,j}^A=0$, then
$$A(i,j)_H\equiv\max\{ k+1 \, |\, h_{i,j}^A=h_{i,j-1}^A=h_{i,j-2}^A=\cdots =h_{i,j-k}^A=0\},
$$
\noindent if $h_{i,j}^A\neq 0$, then we define $A(i,j)_H\equiv 0.$
\end{definition}
\medskip

\begin{definition} Let $\{b_1, b_2 ,\ldots , b_{n-1}\}$ be a set of oriented loops as in Figure \ref{cbase} and let $H=(h_{i,j})=\left(lk(b_{i}, b_{j}^{+})\right)$, $i, j\in\{ 1, 2, \ldots , n-1\}$. If $h_{i,j}^A=0$, then we define 
$$\tilde{A}(i,j)_H\equiv\max\{ k+1 \, |\, h_{i,j}^A=h_{i-1,j}^A=\cdots =h_{i-k,j}^A=0,\, k\in{\mathbb Z}_{\geqq 0}, 0<i-k\},
$$
if $h_{i,j}^A\neq 0$, then we define $\tilde{A}(i,j)_H\equiv 0.$
\end{definition}
\medskip

\begin{figure}
\begin{center}
\includegraphics[clip]{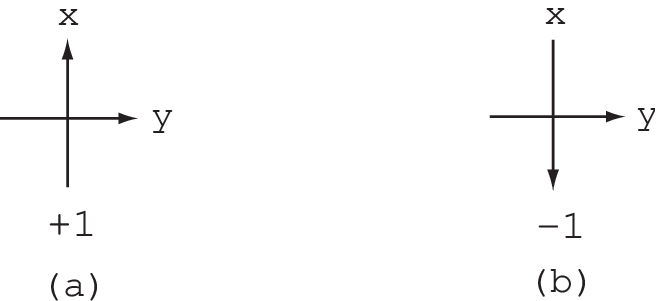}
\end{center}
\caption{}
\label{inter}
\end{figure}

\begin{definition}(cf. ~\cite{7}) Let $S$ be a Seifert surface for a link $L$ and let $x$, $y\in H_1(\textrm{Int}(S))$ be representative loops which intersect transversely. To each intersection of $x$ and $y$ that looks like (a) in Figure \ref{inter} on the positive side of $S$, we assign the value $+1$, whereas to each intersection of $x$ and $y$ that looks like (b) in Figure \ref{inter} on the positive side of $S$, we assign the value $-1$. The intersection number of $x$ and $y$, denoted by $i(x,y)$, is the sum of these signs over all the intersections of $x$ and $y$.
\end{definition}
\medskip

Next we show some properties of a submatrix $H$.

\begin{theorem} \label{Th2} Let $\{ b_1, b_2, \ldots , b_{n-1}\}$ be a set of oriented loops as in the proof of Theorem 1 and let $H=(h_{i,j})=\left(lk(b_{i}, b_{j}^{+})\right)$, $i, j\in\{ 1, 2, \ldots , n-1\}$. Then the following properties hold:

\begin{enumerate}
\item[$(1)$] If $j=i+1$, then $h_{i,j}^A=0$ or $1$, otherwise, $|h_{i,j}^A|=0$ or $1$. 
\item[$(2)$] If $A(i,j)_H$ is even and $i=j-A(i,j)_H$ (resp. $i\neq j-A(i,j)_H$), then $h_{i,j+1}^A=1$ or $0$. (resp. $h_{i,j+1}^A=h_{i,j-A(i,j)_H}^A$ or $h_{i,j+1}^A=0$).
\item[$(3)$] If $A(i,j)_H$ is odd and $i=j-A(i,j)_H$ (resp. $i\neq j-A(i,j)_H$), then $h_{i,j+1}^A=-1$ or $0$ (resp. $h_{i,j+1}^A=-h_{i,j-A(i,j)_H}^A$ or $h_{i,j+1}^A=0$). 
\end{enumerate}
\end{theorem}

\begin{proof} $(1)$ By our choice of loops $\{ b_1, b_2, \ldots , b_{n-1}\}$, we see that $i(b_{i},b_{i+1})$ equals $0$ or $1$ for any $i\in\{1,2,\ldots , n-2\}$. So $h_{i,i+1}^A=0$ or $1$. 

If $j\neq i+1$, then $i(b_{i},b_{j})$ equals $0$ or $1$ or $-1$ for any $i$ and $j$, $i\neq j$, $i, j\in\{1, 2, \ldots , n-1\}$. So $h_{i,j}^A=0$ or $1$ or $-1$. Therefore $(1)$ holds. 

\medskip

\noindent $(2)$ Suppose that $A(i,j)_H$ is even. If $i\neq j-A(i,j)_H$, then we see that $i(b_i, b_{j+1})$ equals $i(b_i, b_{j-A(i,j)_H})$ or $0$. So $h_{i,j+1}^A=h_{i,j-A(i,j)_H}^A$ or $h_{i,j+1}^A=0$. If $i\neq j-A(i,j)_H$, then by $(1)$, we see that $h_{i,j+1}^A=1$ or $0$. Therefore $(2)$ holds.

\medskip

\noindent $(3)$ Suppose that $A(i,j)_H$ is odd. If $i\neq j-A(i,j)_H$, then we see that $i(b_i, b_{j+1})$ equals $-i(b_i, b_{j-A(i,j)_H})$ or $0$. So $h_{i,j+1}^A=-h_{i,j-A(i,j)_H}^A$ or $h_{i,j+1}^A=0$. If $i\neq j-A(i,j)_H$, then by $(1)$, we see that $h_{i,j+1}^A=1$ or $0$. Therefore $(3)$ holds. 
\end{proof}

\begin{theorem} \label{Th3} Let $\{ b_1, b_2, \ldots , b_{n-1}\}$ be a set of oriented loops as in the proof of Theorem 1 and let $H=(h_{i,j})=\left(lk(b_{i}, b_{j}^{+})\right)$, $i, j\in\{ 1, 2, \ldots , n-1\}$. Then the following properties hold:

\begin{enumerate}
\item[$(1)$] If $\tilde{A}(i,j)_H$ is even and $i=j-\tilde{A}(i,j)_H$ (resp. $i\neq j-\tilde{A}(i,j)_H$), then $h_{i,j+1}^A=1$ or $0$. (resp. $h_{i,j+1}^A=h_{i,j-\tilde{A}(i,j)_H}^A$ or $h_{i,j+1}^A=0$).
\item[$(2)$] If $\tilde{A}(i,j)_H$ is odd and $i=j-\tilde{A}(i,j)_H$ (resp. $i\neq j-\tilde{A}(i,j)_H$), then $h_{i,j+1}^A=-1$ or $0$ (resp. $h_{i,j+1}^A=-h_{i,j-\tilde{A}(i,j)_H}^A$ or $h_{i,j+1}^A=0$). 
\end{enumerate}
\end{theorem}

\begin{proof} From condition $(2)$ in Theorem \ref{Th2}, we see that condition $(1)$ of this statement holds. Similarly, From condition $(3)$ in Theorem \ref{Th2}, we see that condition $(2)$ of this statement holds.

\end{proof}

\begin{remark}
By choosing the reverse ordering of the base of a Seifert surface for $L$, we may have the following relation: $h_{i-1,i}^A=0$ or $-1$ for $i=2,3,\ldots , n-1$.\end{remark}
\bigskip
\begin{proposition} Let $\{ b_1, b_2, \ldots , b_{n-1}\}$ be a set of oriented loops as in the proof of Theorem 1 and let $H=(h_{i,j})=(lk(b_{i}, b_{j}^{+}))$, $i, j\in\{ 1, 2, \ldots , n-1\}$. Then the following properties hold:
$$|h_{i,j}^A-h_{i,j+1}^A|=0 \,\,\,\textrm{or} \,\,\, 1,$$ \noindent and
$$|h_{i,j}^A-h_{i-1,j}^A|=0 \,\,\,\textrm{or} \,\,\, 1.$$
\end{proposition}

\begin{proof} 
This is an immediate consequence of Theorem \ref{Th2}.
\end{proof}

\begin{proposition} Let $\{ b_1, b_2, \ldots , b_{n-1}\}$ be a set of oriented loops as in the proof of Theorem 1 and let $H=(h_{i,j})=(lk(b_{i}, b_{j}^{+}))$, $i, j\in\{ 1, 2, \ldots , n-1\}$. If $h_{i,j}^A\neq 0$, then $|\{\, k\in \mathbb{N}\, |\, h_{k,j}^A\neq 0,\,i\leqq k\leqq j-1\}|=|\{\, k\in \mathbb{N}\, |\, h_{i,k}^A\neq 0, \,i+1\leqq k\leqq j\}|\mod 2.$ 
\end{proposition}

\begin{proof} Let $h_{i,j}^A\neq 0$. If $h_{i,j}^A=1$, then, we see that $|\{\, k\in \mathbb{N}\, |\, h_{k,j}^A=0,\,i+1\leqq k\leqq j-1\}|$ and $|\{\, k\in \mathbb{N}\, |\, h_{i,k}^A=0, \,i+1\leqq k\leqq j-1\}|$ are both even. If $h_{i,j}^A=-1$, then, we see that $|\{\, k\in \mathbb{N}\, |\, h_{k,j}^A=0,\,i+1\leqq k\leqq j-1\}|$ and $|\{\, k\in \mathbb{N}\, |\, h_{i,k}^A=0, \,i+1\leqq k\leqq j-1\}|$ are both odd. Hence, $|\{\, k\in \mathbb{N}\, |\, h_{k,j}^A\neq 0,\,i\leqq k\leqq j-1\}|=|\{\, k\in \mathbb{N}\, |\, h_{i,k}^A\neq 0, \,i+1\leqq k\leqq j\}|\mod 2.$  
\end{proof}

\bigskip

\section{A characterization of $2$-component Brunnian links via Seifert matrices}
\medskip

Let $m$ $(\geqq 2)$ be an integer and let $r_{1,m}, r_{2,m}, \ldots , r_{l,m}$ be non-negative integers less than $m$. Then, we consider the following finite sequence of non-negative integers and the set of them, respectively:
$$
( r_{1,m}, r_{2,m}, \ldots , r_{l,m} )_m, 
$$
$$
{\mathcal S}_m\equiv\Big\{ \, ( r_{1,m}, r_{2,m}, \ldots , r_{l,m} )_m\,\, \bigm| \,\,
0\leqq r_{i,m}\leqq m-1,\, i=1,2,\ldots ,l, \,l\in{\mathbb N}\, \Big\}.
$$

\begin{definition}
Let $m$ $(\geqq 2)$ and $j$ $(\leqq m-1)$ be an integer. Let $s_m=( r_{1,m}, r_{2,m}, \ldots , r_{l,m} )_m\in {\mathcal S}_m$. If the following equality holds modulo $2$, then we define $\delta (s_m, j)=m$, otherwise, we define $\delta (s_m, j)=0$:
$$\bigl|\bigl\{\,\, i\,\,|\,\, r_{i,m}=j, \, i\in\{1,2,\ldots , l\}\,\,\bigr\}\bigr|\equiv 1 \quad\mod 2.
$$
\end{definition}
\smallskip

We need some definitions and notations. First define the following set of the sequences $(0)_2$ and $(1)_2$:
$$
{\mathcal S}(2)\equiv \bigl\{ (0)_2, (1)_2 \bigr\} \subset {\mathcal S}_2.
$$ 

\noindent Next take $s_2=( r_{1,2} )_2\in{\mathcal S}(2)$ and define the following set 
$$
{\mathcal S}( s_2, 3)=
\left\{
\begin{array}{@{\,}ll}
 e_{3}(2, v_1, l),\,\, e_{3}(2, v_1, r)\, ; \\
 e_{3}(2, v_2, l),\, e_{3}(1, v_2, r)\, ; \\
 e_{3}(1, v_1, l),\, e_{3}(0, v_1, r)
\end{array}
\right\}
\subset {\mathcal S}_{3}, 
$$
\smallskip

\noindent where 
\begin{eqnarray*}
e_{3}(2, v_1, l) &=& (0)_{3}, \\
e_{3}(2, v_1, r) &=& (2)_{3}, \\ 
e_{3}(2, v_2, l) &=& (r_{1,2}, 2 )_{3}, \\
e_{3}(1, v_2, r) &=& (r_{1,2}, 1)_{3}, \\
e_{3}(1, v_1, l) &=& (r_{1,2}, \delta(s_2, 1), 1 )_{3}, \\
e_{3}(0, v_1, r) &=& (r_{1,2}, \delta(s_2, 1), 0 )_{3}.
\end{eqnarray*}

\medskip

\noindent Similarly, take $s_2=( r_{1,2} )_2\in{\mathcal S}(2)$ and take $s_i=( r_{1,i}, r_{2,i},\ldots , r_{l_i,i} )_i\in{\mathcal S}( s_2,\ldots , s_{i-1}, i)$, $i=3, 4, \ldots , n$, $n\geqq 3$ and define the following set 
$${\mathcal S}( s_2, s_3,\ldots ,s_n, n+1)=
\left\{
\begin{array}{@{\,}ll}
 e_{n+1}(n, v_1, l),\,\, e_{n+1}(n, v_1, r)\, ; \\
 e_{n+1}(n, v_2, l),\, e_{n+1}(n-1, v_2, r)\, ; \\
 e_{n+1}(n-1, v_1, l),\, e_{n+1}(n-2, v_1, r)\, ; \\
 e_{n+1}(n-2, v_2, l),\, e_{n+1}(n-3, v_2, r)\, ; \\
   \quad\quad\quad\quad\vdots\quad\quad\quad\quad \vdots \\
 e_{n+1}(4, v, l),\, e_{n+1}(3, v, r)\, ; \\
 e_{n+1}(3, v', l),\, e_{n+1}(2, v', r)\, ; \\
 e_{n+1}(2, v, l),\,  e_{n+1}(1, v, r)\, ; \\
 e_{n+1}(1, v', l),\, e_{n+1}(0, v', r)
\end{array}
\right\}
\subset {\mathcal S}_{n+1}, 
$$
\smallskip

\noindent where 
\begin{eqnarray*}
e_{n+1}(n, v_1, l) &=& (0)_{n+1}, \\
e_{n+1}(n, v_1, r) &=& (n)_{n+1}, \\ 
e_{n+1}(n, v_2, l) &=& (r_{1,n}, r_{2,n}, \ldots , r_{l_n, n}, n )_{n+1}, \\
e_{n+1}(n-1, v_2, r) &=& (r_{1,n}, r_{2,n}, \ldots , r_{l_n, n}, n-1)_{n+1}, \\
e_{n+1}(n-1, v_1, l) &=& (r_{1,n}, r_{2,n}, \ldots , r_{l_n, n}, r_{1,n-1}, r_{2,n-1}, \ldots , r_{l_{n-1}, n-1}, \delta(s_n, n-1), n-1 )_{n+1}, \\
e_{n+1}(n-2, v_1, r) &=& (r_{1,n}, r_{2,n}, \ldots , r_{l_n, n}, r_{1,n-1}, r_{2,n-1}, \ldots , r_{l_{n-1}, n-1}, \delta(s_n, n-1), n-2 )_{n+1}, \\
e_{n+1}(n-2, v_2, l) &=& (r_{1,n}, r_{2,n}, \ldots , r_{l_n, n}, r_{1,n-1}, r_{2,n-1}, \ldots , r_{l_{n-1}, n-1}, \delta(s_n, n-1), \\
    & & {} r_{1,n-2}, r_{2,n-2}, \ldots, r_{l_{n-2},n-2}, 
                     \delta(s_{n}, n-2), \delta(s_{n-1}, n-2), n-2 )_{n+1}, \\
e_{n+1}(n-3, v_2, r) &=& (r_{1,n}, r_{2,n}, \ldots , r_{l_n, n}, r_{1,n-1}, r_{2,n-1}, \ldots , r_{l_{n-1}, n-1}, \delta(s_n, n-1), \\
    & & {} r_{1,n-2}, r_{2,n-2}, \ldots, r_{l_{n-2},n-2}, \delta(s_{n}, n-2), \delta(s_{n-1}, n-2), n-3 )_{n+1}, \\
\vdots\quad\quad &\vdots& {} \quad\quad\vdots  \\
e_{n+1}(2, v, l) &=& (r_{1,n}, r_{2,n}, \ldots , r_{l_n, n}, r_{1,n-1}, r_{2,n-1}, \ldots , r_{l_{n-1}, n-1}, \delta(s_n, n-1), \\
& &  r_{1,n-2}, r_{2,n-2}, \ldots, r_{l_{n-2},n-2}, \delta(s_{n}, n-2), \delta(s_{n-1}, n-2),\\
& & r_{1,n-3}, r_{2,n-3}, \ldots, r_{l_{n-3},n-3}, \delta(s_{n}, n-3), \delta(s_{n-1}, n-3), \delta(s_{n-2}, n-3), \\
& & {} \quad\quad\quad\quad\vdots \quad\quad \vdots\\
& & r_{1,2}, r_{2,2}, \ldots, r_{l_{2},2}, \delta(s_{n}, 2), \delta(s_{n-1}, 2), \ldots , \delta(s_{4}, 2), \delta(s_{3}, 2), 2) \\
e_{n+1}(1, v, r) &=& (r_{1,n}, r_{2,n}, \ldots , r_{l_n, n}, r_{1,n-1}, r_{2,n-1}, \ldots , r_{l_{n-1}, n-1}, \delta(s_n, n-1), \\
& & r_{1,n-2}, r_{2,n-2}, \ldots, r_{l_{n-2},n-2}, \delta(s_{n}, n-2), \delta(s_{n-1}, n-2), \\
& & r_{1,n-3}, r_{2,n-3}, \ldots, r_{l_{n-3},n-3}, \delta(s_{n}, n-3), \delta(s_{n-1}, n-3), \delta(s_{n-2}, 2), \\
& & {} \quad\quad\quad\quad\vdots \quad\quad \vdots\\
& & r_{1,2}, r_{2,2}, \ldots, r_{l_{2},2}, \delta(s_{n}, 2), 
\delta(s_{n-1}, 2), \ldots , \delta(s_{4}, 2), \delta(s_{3}, 2), 1), \\
e_{n+1}(1, v', l) &=& (r_{1,n}, r_{2,n}, \ldots , r_{l_n, n}, r_{1,n-1}, r_{2,n-1}, \ldots , r_{l_{n-1}, n-1}, \delta(s_n, n-1), \\
& & r_{1,n-2}, r_{2,n-2}, \ldots, r_{l_{n-2},n-2}, \delta(s_{n}, n-2), \delta(s_{n-1}, n-2), \\
& & r_{1,n-3}, r_{2,n-3}, \ldots, r_{l_{n-3},n-3}, \delta(s_{n}, n-3), \delta(s_{n-1}, n-3), \delta(s_{n-2}, n-3), \\
& & {} \quad\quad\quad\quad\vdots \quad\quad \vdots\\
& & r_{1,2}, r_{2,2}, \ldots, r_{l_{2},2}, \delta(s_{n}, 2), 
\delta(s_{n-1}, 2), \ldots , \delta(s_{4}, 2), \delta(s_{3}, 2) \\
& & \delta(s_{n}, 1), \delta(s_{n-1}, 1), \ldots , \delta(s_{3}, 1), \delta(s_{2}, 1), 1) \\
e_{n+1}(0, v', r) &=& (r_{1,n}, r_{2,n}, \ldots , r_{l_n, n}, r_{1,n-1}, r_{2,n-1}, \ldots , r_{l_{n-1}, n-1}, \delta(s_n, n-1), \\
& & r_{1,n-2}, r_{2,n-2}, \ldots, r_{l_{n-2},n-2}, \delta(s_{n}, n-2), \delta(s_{n-1}, n-2), \\
& & r_{1,n-3}, r_{2,n-3}, \ldots, r_{l_{n-3},n-3}, \delta(s_{n}, n-3), \delta(s_{n-1}, n-3), \delta(s_{n-2}, 2), \\
& & {} \quad\quad\quad\quad\vdots \quad\quad \vdots\\
& & r_{1,2}, r_{2,2}, \ldots, r_{l_{2},2}, \delta(s_{n}, 2), 
\delta(s_{n-1}, 2), \ldots , \delta(s_{4}, 2), \delta(s_{3}, 2), \\
& & \delta(s_{n}, 1), \delta(s_{n-1}, 1), \ldots , \delta(s_{3}, 1), \delta(s_{2}, 1), 0), 
\end{eqnarray*}
\noindent and if $n$ is odd, then $v=v_2$ and $v'=v_1$, otherwise, $v=v_1$ and $v'=v_2.$
\bigskip

For example, 
$$
{\mathcal S}\Bigl( (0)_2, 3\Bigr)=\Bigl\{ (0)_3, (2)_3, (0,2)_3, (0,1)_3, (0,0,1)_3, (0,0,0)_3 \Bigr\} \subset {\mathcal S}_3,
$$
because $e_3(2,v_1,l)=(0)_3$, $e_3(2,v_1,r)=(2)_3$, $e_3(2,v_2,l)=(0,2)_3$, $e_3(1,v_2,r)=(0,1)_3$, $e_3(1,v_1,l)=(0,0,1)_3$, and $e_3(1,v_1,r)=(0,0,0)_3$. And for example, 

$$
{\mathcal S}\Bigl( (1)_2, 3\Bigr)=\Bigl\{ (0)_3, (2)_3, (1,2)_3, (1,1)_3, (1,2,1)_3, (1,2,0)_3 \Bigr\} \subset {\mathcal S}_3,
$$
because $e_3(2,v_1,l)=(0)_3$, $e_3(2,v_1,r)=(2)_3$, $e_3(2,v_2,l)=(1,2)_3$, $e_3(1,v_2,r)=(1,1)_3$, $e_3(1,v_1,l)=(1,2,1)_3$, and $e_3(1,v_1,r)=(1,2,0)_3$.
\medskip

\begin{definition}
Let $m$ $(\geqq 2)$ be an integer, and let $s_m=( r_{1,m}, r_{2,m}, \ldots , r_{l_1,m} )_m$, $s'_m=( r'_{1,m}, r'_{2,m}, \ldots , r'_{l_2,m} )_m\in {\mathcal S}( s_2, s_3,\ldots ,s_{m-1}, m)$. Then $s_m$ and $s'_m$ are said to be equivalent in ${\mathcal S}( s_2, s_3,\ldots ,s_{m-1}, m)$, denoted by $s_m\sim_m s'_m$, if the following equality holds modulo $2$ for any $j=1, 2, \ldots , m-1$:
$$
\bigl|\bigl\{\, i\, |\,\, r_{i,m}=j,\, i=1,2,\ldots , l_1\,\bigr\}\bigr|\equiv \bigl|\bigl\{\, i\, |\,\, r'_{i,m}=j,\, i=1,2,\ldots , l_2\,\bigr\}\bigr| \mod 2.$$
\end{definition}

%\begin{definition}
%An element $s_m=( r_{1,m}, r_{2,m}, \ldots , r_{l,m} )_m$ of ${\mathcal S}( g%_2, s_3,\ldots ,s_m, m)$ is said to be representative, if $r_{1,m}<r_{2,m}< \cd%ots <r_{l,m}$. 
%\end{definition}
\medskip

For example, 
$$
{\mathcal S}\Bigl( (0)_2, 3\Bigr)/\sim_3\, = \Bigl\{ \bigl[(0)_3\bigr], \bigl[(2)_3\bigr], \bigl[(0,1)_3\bigr] \Bigr\} \subset {\mathcal S}_3,
$$
$$
{\mathcal S}\Bigl( (1)_2, 3\Bigr)/\sim_3\, =\Bigl\{ \bigl[(0)_3\bigr], \bigl[(2)_3\bigr], \bigl[(1,2)_3\bigr] \Bigr\} \subset {\mathcal S}_3.
$$
\noindent Because in the set ${\mathcal S}\Bigl( (0)_2, 3\Bigr)$, $(0)_3\sim (0,0,0)_3$, $(2)_3\sim (0,2)_3$, and $(0,1)_3\sim (0,0,1)_3$, and in the set ${\mathcal S}\Bigl( (1)_2, 3\Bigr)$, $(0)_3\sim (1,1)_3$, $(2)_3\sim (1,2,1)_3$, and $(1,2)_3\sim (1,2,0)_3.$
\medskip

Let ${\mathcal V}_{n}$ be the set of $(n-1)$ dimensional vectors whose element is $0$ or $1$ or $-1$, $n=2, 3, \ldots$. Let $[ s_{2}]\in{\mathcal S}(2)/\sim_{2}$ and $[ s_{j}]\in{\mathcal S}( s_2, s_3,\ldots ,s_{j-1}, j)/\sim_{j}$, $j=3,4,\ldots $. Then we define functions
$${\mathcal G}_{2}:{\mathcal S}(2)/\sim_{2} \longrightarrow {\mathcal V}_2, 
$$
$${\mathcal G}_{2}\Bigl(\bigl[ s_{2}\bigr]\Bigr)=(x_{1,2}), 
$$
and
$${\mathcal G}_{s_2, s_3, \ldots , s_{j-1},j}:{\mathcal S}( s_2, s_3,\ldots ,s_{j-1}, j)/\sim_{j} \longrightarrow {\mathcal V}_j, 
$$ 
$${\mathcal G}_{s_2, s_3, \ldots , s_{j-1},j}\Bigl(\bigl[ s_{j}\bigr]\Bigr)=(x_{1,j}, x_{2,j},\ldots , x_{j-1,j}), 
$$

\noindent recursively, as follows: first, let 
$${\mathcal G}_2\Bigl(\bigl[ ( 0 )_2\bigr]\Bigr)=(x_{1,2})=(0), 
$$
$${\mathcal G}_2\Bigl(\bigl[ ( 1 )_2\bigr]\Bigr)=(x_{1,2})=(1).
$$
Next, given $s_2, s_3, \ldots , s_{n-1}$, ${\mathcal G}_{2}$ and ${\mathcal G}_{s_2,s_3,\ldots , s_{j-1}, j}, j=3,4,\ldots , n$, then $x_{k, n+1}=1$ is decided as follows: If $\delta(s_{n+1},k)\neq 0$ and $x_{k, n}=1$ (resp. $-1$), then we define $x_{k, n+1}=1$ (resp. $-1$). If $\delta (s_{n+1},k)\neq 0$, $x_{k, n}=0$, and the following integer 
$$
B(k,n+1)\equiv\max\{l+1|x_{k, n}=x_{k, n-1}=\cdots =x_{k, n-l}=0\}
$$
\noindent is even (resp. odd), then we define $x_{k, n+1}=x_{k, n+1-B(k, n+1)}$ (resp. $-x_{k, n+1-B(k, n+1)}$). If $\delta(s_{n+1},k)=0$, then we define $x_{k, n+1}=0$. Here, $k\in\{1,2,\ldots , n-1\}$. 

If $\delta(s_{n+1}, n)\neq 0$, then we define $x_{n, n+1}=1$. If $\delta(s_{n+1}, n)=0$, we define $x_{n, n+1}=0$.

\medskip

Let $\{ b_1, b_2, \ldots , b_{n-1}\}$ be a set of oriented loops as in the proof of Theorem 1 and let $H=(h_{i,j})=\bigl(lk(b_{i}, b_{j}^{+})\bigr)$, $i, j\in\{ 1, 2, \ldots , n-1\}$. Then $\bigl(h_{1, k}^A, h_{2,k}^A, \ldots , h_{k-1,k}^A\bigr)\in {\mathcal V}_k$ for each $k\in\{2, 3, \ldots , n-1\}$. 
\medskip

In the proof of Theorem 1, the set $D_2\cap b_i$ is an oriented arc. The initial point of the oriented arc $D_2\cap b_i$ is denoted by $I_i$ and the terminal point of the oriented arc $D_2\cap b_i$ is denoted by $T_i$. Also, the initial point of the oriented arc $D^2_2\cap b_i$ is denoted by $I_i$ and the terminal point of the oriented arc $D^2_2\cap b_i$ is denoted by $T_i$.

\begin{theorem}\label{necess} Any Seifert matrix of a 2-component Brunnian link is S-equivalent to a matrix $M_L$ of the following form:
$$
M_L=\left(
\begin{array}{c|c}
E & F_{n}  \\ \hline
0 & H \\
\end{array}
\right),
$$
\noindent and there are $s_2\in{\mathcal S}(2)$ and $s_j\in{\mathcal S}(s_2, s_3, \ldots , s_{j-1}, j)$, $j=3, 4, \ldots , n-1$ such that 
$${\mathcal G}_2\bigl([s_2]\bigr)=\bigl(h_{1,2}^A\bigr), 
$$ and 
$${\mathcal G}_{s_2, s_3, \ldots , s_{j-2}, j-1}\bigl([s_{j-1}]\bigr)=\bigl( h_{1,j-1}^A, h_{2,j-1}^A, \ldots , h_{j-2,j-1}^A\bigr).$$

\end{theorem}

\begin{proof} By $(1)$ in Theorem \ref{Th2}, $h_{1,2}^A=0$ or $1$. By the definition of $\mathcal{G}_2$, there is an element $s_2$ in $\mathcal{S}(2)$ such that $\mathcal{G}_2([s_2])=(h_{1,2}^A)$, i.e., if $h_{1,2}^A=0$, then by taking the element $(0)_2$ of ${\mathcal S}(2)$ as $s_2$, we obtain ${\mathcal G}_2([s_2])=0=h_{1,2}^A$, otherwise, by taking the element $(1)_2$ of ${\mathcal S}(2)$ as $s_2$, we obtain ${\mathcal G}_2([s_2])=1=h_{1,2}^A$. 
\medskip

Now, a pair of arcs $b_{i-1}$ and $b_i$ is arranged on the positive side of $D_2^2$ as in Figure \ref{Th4-3}, where $i=2,3,\ldots , n-1$. The singularity of $D_1^2$ and $D_2^2$ which $b_{i-1}$ and $b_i$ go through is denoted by $c_{i-1}$, where $i=2,3,\ldots , n-1$ (see Figure \ref{Th4-3}). The singularity of $D_1^2$ and $D_2^2$ which only $b_1$ (resp. $b_{n-1}$) goes through is denoted by $c_1$ (resp. $c_{n-1}$). 
\medskip

Given $b_1$ and $b_2$ on $D_2$, we construct $b_3$ from the information of $b_1$ and $b_2$. On $c_2$, $I_3$ (or $T_3$) is adjacient clockwisely to $I_2$ (or $T_2$). So we translate $T_3$ (or $I_3$) starting from a neighbourhood of $I_3$ (or $T_3$) along $b_1$, $b_2$ or singularities to fix its position on $\partial D_2^2$ except $c_1$ and $c_2$. Then from our choice of bases we see that $b_3\cap D_2$ is one of three arcs obtained by the following three ways $(1)-(3)$:
\medskip

(1) Translate $T_3$ (or $I_3$) clockwisely a little along $\partial D_2$ starting from a neighbourhood of $I_3$ (or $T_3$) to fix its position on $\partial D_2$ which is adjacient to $T_3$ (or $I_3$). Then $T_3$ (or $I_3$) does not transverse any arc $b_i$, $i=1,2$, while translating it from the starting point to the fixed point. We write arcs $D_2\cap b_j$, $j=1,2$ as a sequence $(j_1, j_2, \ldots )_3$ in turns that $T_3$ transverses while translating it from the starting point to the fixed point. Then the arc $b_3\cap D_2^2$ does not intersect any arc $b_i$, $i=1,2$. So we write $(0)_3$. 
\medskip

(2) Translate $T_3$ (or $I_3$) counterclockwisely starting from a neighbourhood of $I_3$ (or $T_3$) through a neighbourhood of $I_2$ (or $T_2$) and translate it along $b_2$ until it arrive at $T_2$ (or $I_2$). After translating $T_3$ (or $I_3$) a little along $\partial D_2$ clockwisely from $T_2$ (or $I_2$) to fix $T_3$ (or $I_3$). Then $T_3$ (or $I_3$) goes once through $b_i$ which intersect $b_2$, $i=1$ and goes once through $b_2$. Let $s_2=(r_{1,2})_2$. As a result, if $r_{1,2}=0$, then $b_3\cap D_2$ intersects $b_2$ once, and if $r_{1,2}\neq 0$, i.e., $r_{1,2}=1$, then $b_3\cap D_2$ intersects $b_1$ once and intersects $b_2$ once. We denote it by $(r_{1,2}, 2)_3$. 
\medskip

(3) Similarly, translate $T_3$ (or $I_3$) counterclockwisely starting from a neighbourhood of $I_3$ (or $T_3$) to goes through a neighbourhood of $I_2$ (or $T_2$) and translate it along $b_2$ until it arrive at $T_2$ (or $I_2$). And proceed to translate it along $c_1$ until it arrives at $T_1$ (or $I_1$) and proceed to translate along $b_1$ until it arrives at $I_1$ (or $T_1$). After translating $T_3$ (or $I_3$) a little along $\partial D_2$ clockwisely from $T_1$ (or $I_1$) to fix $T_3$ (or $I_3$). Then $T_3$ (or $I_3$) goes once through $b_i$ which intersect $b_2$, $i=1$ and goes once through $b_i$ which intersect $b_1$, $i=2$ and goes through $b_1$. Then we denote it by $(r_{1,2}, \sigma(s_2,1), 1)_3$. 
\medskip

If $b_3$ is realized by the way (1), then we set $s_3=(0)_3$. Then $s_3\in{\mathcal S}(s_2, 3)$ and $\mathcal{G}_{s_2,3}([s_3])=(0,0)=(h_{1,3}^A, h_{2,3}^A)$ by the definition of $\mathcal{G}_{s_2, 3}$. If $b_3$ is realized by the way (2), then we set $s_3=(r_{1,2}, 2)_3$. Then $s_3\in{\mathcal S}(s_2, 3)$ and $\mathcal{G}_{s_2,3}([s_3])=(h_{1,3}^A, h_{2,3}^A)$ by the definition of $\mathcal{G}_{s_2, 3}$. If $b_3$ is realized by the way (3), then we set $s_3=(r_{1,2}, \sigma(s_2,1), 1)_3$. Then $s_3\in{\mathcal S}(s_2, 3)$ and $\mathcal{G}_{s_2,3}([s_3])=(h_{1,3}^A, h_{2,3}^A)$ by the definition of $\mathcal{G}_{s_2, 3}$. 
\medskip

By the obove way, inductively, we construct $b_i$, $i=4,5,\ldots , n-1$. And we set $s_4\in{\mathcal S}(s_2, s_3, 4)$, $s_5\in{\mathcal S}(s_2, s_3, s_4,5)$, $\ldots$, $s_{n-1}\in{\mathcal S}(s_2, s_3,\ldots, s_{n-2}, n-1)$. As a result, by its construction, we see that for any $j=3,4,\ldots , n-1$,
$$
{\mathcal G}_{s_2, s_3, \ldots , s_{j-2}, j-1}\bigl([s_{j-1}]\bigr)=\bigl( h_{1,j-1}^A, h_{2,j-1}^A, \ldots , h_{j-2,j-1}^A\bigr).
$$
\end{proof}

\begin{proposition} Let $H=(h_{i,j})=\left(lk(b_{i}, b_{j}^{+})\right)$, $i, j\in\{ 1, 2, \ldots , n-1\}$. Two end points $T_i$ and $T_j$ are adjacent in $\partial D_2$ if and only if $e_n(i, v, w_1)\sim_{n}e_n(j, v, w_2)$. Two end points $I_i$ and $I_j$ are adjacent in $\partial D_2$ if and only if $e_n(i, v', w_1)\sim_{n}e_n(j, v', w_2)$. Two end points $I_i$ and $T_j$ are adjacent in $\partial D_2$ if and only if $e_n(i, v, w_1)\sim_{n}e_n(j, v', w_2)$, where $v\neq v'$, $v, v'\in\{v_1, v_2\}$ and $w_1\neq w_2$, $w_1, w_2\in\{l, r\}.$
\end{proposition}

\begin{proof} This is an immediate consequence of the diagram of $D_2$.

\end{proof}

\begin{theorem}\label{suff}
Let $M$ be a square matrix having the following form: 
$$
M=\left(
\begin{array}{c|c}
E & F_{n}  \\ \hline
0 & H \\
\end{array}
\right).
$$
If the following equalities hold for some $s_2\in{\mathcal S}(2)$ and $s_j\in{\mathcal S}(s_2, s_3, \ldots , s_{j-1}, j)$, $j=3, 4, \ldots , n-1$, then there are a 2-component Brunnian link $L$ and a Seifert matrix $M_L$ for $L$ such that $M_L=M$: 
$${\mathcal G}_2\bigl([s_2]\bigr)=\bigl(h_{1,2}^A\bigr), 
$$ 
and 
$${\mathcal G}_{s_2, s_3, \ldots , s_{j-2}, j-1}\bigl([s_{j-1}]\bigr)=\bigl( h_{1,j-1}^A, h_{2,j-1}^A, \ldots , h_{j-2,j-1}^A\bigr).$$ 
\end{theorem}

\begin{proof}
First, we consider the case where $n=1$. Then, matrices satisfying the assumption in Theorem 4 are the two matrices $(1)$ and $(-1)$. It is easy to see that there is a Seifert surface whose boundary is a 2-component Brunnian link that realizes each matrix $(1)$ or $(-1)$.

Next, we consider the case where $n=2$. Then, matrices satisfying the assumption in Theorem 4 are the following:
$$
\left(
\begin{array}{ccc}
1 & 0 & 1 \\ 
0 & 1 & 1 \\ 
0 & 0 & h_{1,1} \\
\end{array}
\right), 
\left(
\begin{array}{ccc}
1 & 0 & 1 \\ 
0 & -1 & 1 \\ 
0 & 0 & h_{1,1} \\
\end{array}
\right), 
\left(
\begin{array}{ccc}
-1 & 0 & 1 \\ 
0 & -1 & 1 \\ 
0 & 0 & h_{1,1} \\
\end{array}
\right).
$$
\noindent We denote the above matrices, by $M_1$, $M_2$, and $M_3$ from the left. It is easy to see that for any integer $h_{1,1}$, there exists a Seifert surface whose boundary is a 2-component Brunnian link that realizes $M_1$.  Namely, for each integer $h_{1,1}$, a 2-component Brunnian link whose Seifert surface is $M_1$ is the link having $h_{1,1}+1$ twists. Similarly, for each integer $h_{1,1}\neq 0$, a 2-component Brunnian link whose Seifert surface is $M_2$ is the link having $h_{1,1}$ twists and for each integer $h_{1,1}$, a 2-component Brunnian link whose Seifert surface is $M_3$ is the link having $h_{1,1}-1$ twists.

Also, when $h_{1,1}\neq 0$, there is a Seifert surface that realizes $M_2$ (cf. Figure ). Note that although the Alexander-Conway polynomials of the links whose Seifert surfaces are as shown in Figure , say $L_1$ and $L_2$, are equal to zero, $L_1$ and $L_2$ are both non-trivial and hence they are Brunnian. Because the Jones polynomials of $L_1$ and $L_2$ are non-trivial as follows:
\begin{eqnarray*}
V_{L_1}(t) &=& t^{-15/2}-2t^{-13/2}+2t^{-11/2}-2t^{-9/2}+2t^{-7/2}-6t^{-5/2}+
9t^{-3/2} \\
& & -10t^{-1/2}+6t^{1/2}-3t^{3/2}+3t^{5/2}-7t^{7/2}+t^{9/2}, \\
V_{L_2}(t) &=& -t^{-9/2}+t^{-7/2}+2t^{-5/2}-3t^{-3/2}+5t^{-1/2}-8t^{1/2}+10t^{3/2} \\
& & -5t^{5/2}+2t^{7/2}-2t^{9/2}+2t^{11/2}-2t^{13/2}+t^{15/2}.
\end{eqnarray*}

\medskip

Next, we consider the case where $n\geqq 3$. For any matrix $M$ satisfying the assumption in Theorem 4, we construct a 2-component Brunnian link $L$ such that $M_L=M$. Let $M$ be a matrix satisfying the assumption of Theorem 4. 
\smallskip

Let $E=(e_{i,j})$, $i$, $j\in\{1,2,\ldots , n\}$. If $e_{i,i}=e_{i+1,i+1}=1$, then choose arcs $b_i'$, $b_{i+1}'$, $c_i'$, and $c_{i+1}'$ on a disk, say $D_1^2$ as the top and left side diagram in Figure . If $e_{i,i}=e_{i+1,i+1}=-1$, then choose arcs $b_i'$, $b_{i+1}'$, $c_i'$, and $c_{i+1}'$ on $D_1^2$ as the top and right side diagram in Figure . If $e_{i,i}=1$ and $e_{i+1,i+1}=-1$, then choose arcs $b_i'$, $b_{i+1}'$, $c_i'$, and $c_{i+1}'$ on $D_1^2$ as the bottom and left side diagram in Figure . If $e_{i,i}=-1$ and $e_{i+1,i+1}=1$, then choose arcs $b_i'$, $b_{i+1}'$, $c_i'$, and $c_{i+1}'$ on $D_1^2$ as the bottom and right side diagram in Figure . Thus, we choose an arc $b_i'$ for any $i\in\{1,2,\ldots , n-1\}$ whose end points are on $Int(c_i')$ and $Int(c_{i+1}')$, where $Int(*)$ is the set of interior points of $*$. As orientations of $b_1'$ and $b_2'$, we may take the orientations as in Figure . So the orientations of $b_3'$, $b_4', \ldots , b_{n-1}'$ are all decided so that for any $i\in\{1,2,\ldots , n-1\}$, the orientations of $b'_i$ and $b'_{i+1}$ are the same as one (the left side or the right side diagram) of the diagrams of Figure \ref{connect2}. 
\smallskip

Now by the assumption, the following equalities hold for some $s_2\in{\mathcal S}(2)$ and $s_j\in{\mathcal S}(s_2, s_3, \ldots , s_{j-1}, j)$, $j=3, 4, \ldots , n-1$: 
$${\mathcal G}_2\bigl([s_2]\bigr)=\bigl(h_{1,2}^A\bigr), 
$$ 
and 
$${\mathcal G}_{s_2, s_3, \ldots , s_{j-2}, j-1}\bigl([s_{j-1}]\bigr)=\bigl( h_{1,j-1}^A, h_{2,j-1}^A, \ldots , h_{j-2,j-1}^A\bigr).$$ 

Therefore, there are oriented arcs $b_1, b_2, \ldots , b_{n-1}$ on a disk, say $D_2$ such that $i(b_i, b_j)=h_{i,j}^A$, where $H=(h_{i,j})$, $i$, $j\in\{1,2,\ldots , n-1\}$. And then deform $D_2$ into a disk, say $D_2^2$ according the elements of $E$ as in Figure \ref{Th4-3}. For example, if $e_{i,i}-1$, then we deform $D_2^2$. We denote the set of the double points by $c_{i}$ (see Figure \ref{Th4-3}). 

\begin{figure}
\begin{center}
\includegraphics[clip]{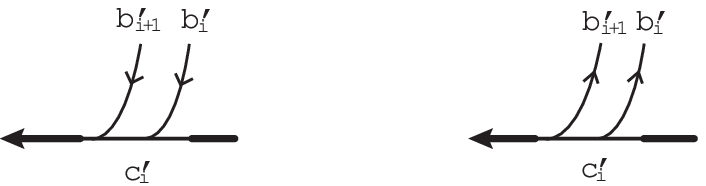}
\end{center}
\caption{}
\label{connect2}
\end{figure}

\begin{figure}
\begin{center}
\includegraphics[clip]{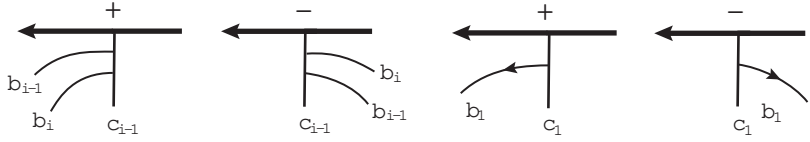}
\end{center}
\caption{arcs on $D_2^2$}
\label{Th4-3}
\end{figure}
\smallskip

Next, we construct a C-complex from the two disks $D_1^2$ and $D_2^2$. First, after tying $D_1^2$ in $\mathbb{S}^3$ by making a clasp singularity that is just $c'_2=c_2$ such that $T'_1=I_1$ according to the element $e_{1,1}$ of $E$ as in Figure . Next, after tying $D_1^2$ and twisting $D_1^2$ ($h_{1,1}+6$ twists), make a clasp singularity that is just $c'_2=c_2$ such that $I'_1=T_1$ and $I'_2=T_2$ according to the element $e_{2,2}$ of $E$ as in Figure . 

Next intertwine a part of $D_1^2$ containing $c_3$ around $b_1\cup b'_1$, $h_{1,2}$ times. And then after tying $D_1^2$ and twisting $D_1^2$ ($h_{2,2}+6-h_{1,2}$ twists), we make a clasp singularity by connecting arcs $c'_3$ and $c_3$ as in Figure . 

Similarly, intertwine a part of $D_1^2$ containing $c_4$ around $b_1\cup b'_1$, $h_{1,3}$ times and intertwine a part of $D_1^2$ containing $c_4$ around $b_2\cup b'_2$, $h_{2,3}$ times. And then after tying $D_1^2$ and twisting $D_1^2$ ($h_{3,3}+6-(h_{1,3}+h_{2,3})$ twists), we make a clasp singularity by connecting arcs $c'_4$ and $c_4$ as in Figure . 

Similarly, we make clasp singularities $c_i=c'_i$ in order from $i=5$ to $i=n$ as follows: intertwine a part of $D_1^2$ containing $c_i$ around $b_j\cup b'_j$, $h_{j,i-1}$ times for each $j=1,2,\ldots , i-1$. And then after tying $D_1^2$ and twisting $D_1^2$ ($h_{i-1,i-1}+6-\sum_{j=1}^{i-1} h_{j,i-1}$ twists), we make a clasp singularity by connecting arcs $c'_i$ and $c_i$ as in Figure . 

In this way clasp singularities $c_1=c'_1$, $c_2=c'_2, \ldots , c_n=c'_n$ are made and we obtain a C-complex, say $C_M$, as in Figure .  

Then the boundary $\partial C_M$ of $C_M$ is a 2-componet Brunnian link. This completes the proof.
\end{proof}

%\begin{corollary} Let $M_L$ be a Seifert matrix for a 2-component Brunnian link% $L$ as in Theorem 1. Then the following equality holds:
%$$
%M_L-tM_L^T=\left(
%\begin{array}{c|c}
%(1-t)E & F_n  \\ \hline

%-tF_{n}^T & H-tH^T \\
%\end{array}
%\right).
%$$
%\noindent Here $H$ is an $(n-1)\times (n-1)$ matrix as in Theorem 1 and $X^T$ d%enotes the transposed matrix of $X$.
%\end{corollary}
%\begin{proof} This is an immediate consequence of Theorem 1.
%\end{proof}

%\begin{remark}
%Let $H^*\equiv H-tH^T=(h_{i,j}^*)$. Then the following conditions hold: if $i=j%$, then $h_{i,j}^*=h_{i,j}(1-t)$; otherwise, $h_{i,j}^*=h_{i,j}(1-t)$ or $h_{i,%j}(1-t)\pm t=h_{j,i}(1-t)\pm 1$.  
%\end{remark}

%\begin{remark} If a 2-component Brunnian link $L$ is algebraically split, then %the submatrix $E^*\equiv (1-t)E=(e_{i,j}^*)$ of $M_L-tM_L^T$ for $L$ as in Coro%llary 2 satisfies the equality: $\sum_{i=1}^n e_{ii}^*=0$.
%\end{remark}

\begin{example} Let $\{ b_1, b_2\}$ be a set of oriented loops as in the proof of Theorem 1 and let $H=(h_{i,j})=(lk(b_{i}, b_{j}^{+}))$, $i, j\in\{ 1, 2\}$. Then $h_{i,j}^A=0$ or $1$.
\end{example}

\begin{example} Let $\{ b_1, b_2, b_3\}$ be a set of oriented loops as in the proof of Theorem 1 and let $H=(h_{i,j})=(lk(b_{i}, b_{j}^{+}))$, $i, j\in\{ 1, 2, 3\}$. Then $(h_{1,2}^A, h_{1,3}^A, h_{2,3}^A)$ is one of the following:
\smallskip

\begin{enumerate}
\item[1.] $(h_{1,2}^A, h_{1,3}^A, h_{2,3}^A)=(0,0,0),$ \\
\item[2.] $(h_{1,2}^A, h_{1,3}^A, h_{2,3}^A)=(0,0,1),$ \\
\item[3.] $(h_{1,2}^A, h_{1,3}^A, h_{2,3}^A)=(0,-1,0),$ \\
\item[4.] $(h_{1,2}^A, h_{1,3}^A, h_{2,3}^A)=(1,0,0),$ \\
\item[5.] $(h_{1,2}^A, h_{1,3}^A, h_{2,3}^A)=(1,1,1),$ \\
\item[6.] $(h_{1,2}^A, h_{1,3}^A, h_{2,3}^A)=(1,0,1).$ \\
\end{enumerate}
\end{example}

\begin{example} Let $\{ b_1, b_2, b_3, b_4\}$ be a set of oriented loops as in the proof of Theorem 1 and let $H=(h_{i,j})=(lk(b_{i}, b_{j}^{+}))$, $i, j\in\{ 1, 2, 3, 4\}$. Then $(h_{1,2}^A, h_{1,3}^A, h_{2,3}^A, h_{1,4}^A, h_{2,4}^A, h_{3,4}^A)$ is one of the following:
\smallskip

\begin{enumerate}
\item[1.] $(h_{1,2}^A, h_{1,3}^A, h_{2,3}^A, h_{1,4}^A, h_{2,4}^A, h_{3,4}^A)=(0,0,0,0,0,0),$ \\
\item[2.] $(h_{1,2}^A, h_{1,3}^A, h_{2,3}^A, h_{1,4}^A, h_{2,4}^A, h_{3,4}^A)=(0,0,0,0,0,1),$ \\
\item[3.] $(h_{1,2}^A, h_{1,3}^A, h_{2,3}^A, h_{1,4}^A, h_{2,4}^A, h_{3,4}^A)=(0,0,0,0,-1,0),$ \\
\item[4.] $(h_{1,2}^A, h_{1,3}^A, h_{2,3}^A, h_{1,4}^A, h_{2,4}^A, h_{3,4}^A)=(0,0,0,1,0,0),$ \\
\item[5.] $(h_{1,2}^A, h_{1,3}^A, h_{2,3}^A, h_{1,4}^A, h_{2,4}^A, h_{3,4}^A)=(0,0,1,0,0,0),$ \\
\item[6.] $(h_{1,2}^A, h_{1,3}^A, h_{2,3}^A, h_{1,4}^A, h_{2,4}^A, h_{3,4}^A)=(0,0,1,0,0,1),$ \\
\item[7.] $(h_{1,2}^A, h_{1,3}^A, h_{2,3}^A, h_{1,4}^A, h_{2,4}^A, h_{3,4}^A)=(0,0,1,0,1,1),$ \\
\item[8.] $(h_{1,2}^A, h_{1,3}^A, h_{2,3}^A, h_{1,4}^A, h_{2,4}^A, h_{3,4}^A)=(0,0,1,1,1,1),$ \\
\item[9.] $(h_{1,2}^A, h_{1,3}^A, h_{2,3}^A, h_{1,4}^A, h_{2,4}^A, h_{3,4}^A)=(0,-1,0,0,0,0),$ \\ 
\item[10.] $(h_{1,2}^A, h_{1,3}^A, h_{2,3}^A, h_{1,4}^A, h_{2,4}^A, h_{3,4}^A)=(0,-1,0,0,0,1),$ \\ 
\item[11.] $(h_{1,2}^A, h_{1,3}^A, h_{2,3}^A, h_{1,4}^A, h_{2,4}^A, h_{3,4}^A)=(0,-1,0,-1,0,1),$ \\
\item[12.] $(h_{1,2}^A, h_{1,3}^A, h_{2,3}^A, h_{1,4}^A, h_{2,4}^A, h_{3,4}^A)=(0,-1,0,-1,-1,0),$ \\
\item[13.] $(h_{1,2}^A, h_{1,3}^A, h_{2,3}^A, h_{1,4}^A, h_{2,4}^A, h_{3,4}^A)=(1,0,0,0,0,0),$ \\
\item[14.] $(h_{1,2}^A, h_{1,3}^A, h_{2,3}^A, h_{1,4}^A, h_{2,4}^A, h_{3,4}^A)=(1,0,0,0,0,1),$ \\
\item[15.] $(h_{1,2}^A, h_{1,3}^A, h_{2,3}^A, h_{1,4}^A, h_{2,4}^A, h_{3,4}^A)=(1,0,0,0,-1,0),$ \\
\item[16.] $(h_{1,2}^A, h_{1,3}^A, h_{2,3}^A, h_{1,4}^A, h_{2,4}^A, h_{3,4}^A)=(1,0,0,-1,-1,0),$ \\
\item[17.] $(h_{1,2}^A, h_{1,3}^A, h_{2,3}^A, h_{1,4}^A, h_{2,4}^A, h_{3,4}^A)=(1,1,1,0,0,0),$ \\
\item[18.] $(h_{1,2}^A, h_{1,3}^A, h_{2,3}^A, h_{1,4}^A, h_{2,4}^A, h_{3,4}^A)=(1,1,1,0,0,1),$ \\
\item[19.] $(h_{1,2}^A, h_{1,3}^A, h_{2,3}^A, h_{1,4}^A, h_{2,4}^A, h_{3,4}^A)=(1,1,1,1,1,1),$ \\
\item[20.] $(h_{1,2}^A, h_{1,3}^A, h_{2,3}^A, h_{1,4}^A, h_{2,4}^A, h_{3,4}^A)=(1,1,1,1,0,0),$ \\
\item[21.] $(h_{1,2}^A, h_{1,3}^A, h_{2,3}^A, h_{1,4}^A, h_{2,4}^A, h_{3,4}^A)=(1,0,1,0,0,0),$ \\
\item[22.] $(h_{1,2}^A, h_{1,3}^A, h_{2,3}^A, h_{1,4}^A, h_{2,4}^A, h_{3,4}^A)=(1,0,1,0,0,1),$ \\
\item[23.] $(h_{1,2}^A, h_{1,3}^A, h_{2,3}^A, h_{1,4}^A, h_{2,4}^A, h_{3,4}^A)=(1,0,1,0,1,1),$ \\
\item[24.] $(h_{1,2}^A, h_{1,3}^A, h_{2,3}^A, h_{1,4}^A, h_{2,4}^A, h_{3,4}^A)=(1,0,1,-1,0,1).$ 
\end{enumerate}
\end{example}

\bigskip

\section{Seifert matrices of $n(\geqq 3)$-component Brunnian links}
\medskip

We denote the following $n\times n$ diagonal matrix denoted by $E_{n}=(e_{i,j})$, $i,j\in\{1 ,\ldots ,n\}$: if $i=j=$odd, then $e_{i,j}=1$, if $i=j=2k$, then $e_{i,j}=-1$ for some $k\in\{1 ,\ldots ,{n\over 2}\}$.

 In the next theorem, we introduce our choice of bases of the first homology group on a Seifert surface for an $n$-component Brunnian link $L$, $n\geqq 3$, and construct a Seifert matrix of it.
 
\begin{theorem} \label{ncbm} For any $n$-component Brunnian link $L$, $n\geqq 3$, there is a Seifert matrix $M_L$ for $L$ of the following form:
$$
M_L=\left(
\begin{array}{c|c|c|c|c|c|c|c}
E_{n_1} & 0 & 0 & 0 & F_{n_1} & 0 & 0 & 0 \\ \hline
0 & E_{n_2} & 0 & 0 & 0 & F_{n_2} & 0 & 0 \\ \hline
0 & 0 & \ddots & 0 & 0 & 0 & \ddots  & 0 \\ \hline
0 & 0 & 0 &  E_{n_{n-1}} & 0 & 0 & 0 &  F_{n_{n-1}} \\ \hline
0 & 0 & 0 & 0 & H_1 & P_{1,2} & \cdots & P_{1,{n-1}} \\ \hline
0 & 0 & 0 & 0 & P_{2,1} & \ddots & \ddots & \vdots \\ \hline 
0 & 0 & 0 & 0 & \vdots & \ddots & H_{n-2} & P_{{n-2},{n-1}} \\ \hline 
0 & 0 & 0 & 0 & P_{{n-1},1} & \cdots & P_{{n-1},{n-2}} & H_{n-1} \\
\end{array}
\right).
$$

\noindent Here, $n_i\in 2\mathbb{Z}_{+}$, and each $H_l$ is a square matrix, $1\leqq i,\, l\leqq n-1$.
\end{theorem}

\begin{proof} Let $L=\cup_{i=1}^n K_i$, and let $C_L=\cup _{i=1}^n D^2_i$ and $S_L$ be a C-complex and a Seifert surface for $L$ obtained according to \lq\lq Algorithm\rq\rq\, in section 3. As a set of oriented loops representing a base of $H_1(\textrm{Int}(S_L))$, we may take 
$${\mathcal L}\equiv\{\,\, a_{l,1}, a_{l,2}, \ldots, a_{l,n_l}, b_{l,1}, b_{l,2},\ldots, b_{l,n_l-1}\,\, | \,\, l=1,2,\ldots, n-1 \,\,\}
$$ 
\noindent such that each $D_l\cap (a_{l,k}\cap b_{m,j})$ is the empty set for $l\neq m$, and each $D_n\cap (b_{l,i}\cap b_{m,j})$ consists of at most one point, and
$$
D_l\cap (b_{l,i}\cap b_{m,j})=
\left\{
\begin{array}{@{\,}ll}

\textrm{ at most one point } & \textrm{if} \,\, \, l=m, \\

\, \, \emptyset & \textrm{if}  \,\, \, l\neq m,
\end{array}
\right.
$$
$$
D_n\cap (a_{l,k}\cap b_{m,j})=
\left\{
\begin{array}{@{\,}ll}

\textrm{ at most one point } & \textrm{if} \,\, \, l=m, \\

\, \, \emptyset & \textrm{if}  \,\, \, l\neq m.
\end{array}
\right.
$$
\noindent Furthermore, 
$$
a_{l,k}\cap b_{l,i}=\left\{
\begin{array}{@{\,}ll}
\textrm{ one point } & \textrm{if} \,\, \, k=i, \, i-1,  \\

\, \, \emptyset & \textrm{if}  \,\, \, k\neq i, \, i-1, \, \textrm{and}\,;
\end{array}
\right.
$$
\noindent $lk(a_{l,k}, b_{l,i}^+)=1$ for any $k=i$ or $i-1$, where, an integer $n_l$ is the number of clasp singularities on the disk $D^2_l$ of $C_L$ and $l,m\in\{1, 2, \ldots , n-1\}$, $k\in\{1, 2, \ldots , n_l\}$, $j\in\{1, 2, \ldots , n_m-1\}$, and $i\in\{1, 2, \ldots , n_l-1\}$.  

By using the base ${\mathcal L}$ of $H_1(\textrm{Int}(S_L))$, we obtain a Seifert matrix $M=(m_{i,j})$ for $L$, $i,j\in\{1, 2, \ldots, 2n-1\}$ as follows: 
$$
m_{i,j}=\left\{
\begin{array}{@{\,}ll}

lk(a_{I_1}, a_{J_1}^{+}) & \, \textrm{if} \,\,\,  \hspace{0.2cm}

\sum_{i=0}^k n_i+1\leqq i\leqq\sum_{i=0}^{k+1} n_i\,\, \textrm{and}  \\

& \hspace{1.5cm} \sum_{i=0}^l n_i+1\leqq j\leqq\sum_{i=0}^{l+1} n_i,  \\

& \\

lk(a_{I_2}, b_{J_2}^{+}) & \, \textrm{if} \,\,\, \hspace{0.2cm} \sum_{i=0}^k n_i+1\leqq i\leqq\sum_{i=0}^{k+1} n_i\,\, \textrm{and} \\

& \hspace{0.7cm} \sum_{i=0}^{l}(n_i-1)+1\leqq j-\sum_{i=0}^{n-1} n_i\leqq\sum_{i=0}^{l+1}(n_i-1),  \\

& \\

lk(b_{I_3}, a_{J_3}^{+}) & \, \textrm{if} \,\,\, \hspace{0.3cm} \sum_{i=0}^{k}(n_i-1)+1\leqq i-\sum_{i=0}^{n-1} n_i\leqq\sum_{i=0}^{k+1}(n_i-1)\,\, \textrm{and} \\

& \hspace{1.5cm} \sum_{i=0}^l n_i+1\leqq j\leqq\sum_{i=0}^{l+1} n_i, \\

& \\

lk(b_{I_4}, b_{J_4}^{+}) & \, \textrm{if} \,\,\, \hspace{0.3cm} \sum_{i=0}^{k}(n_i-1)+1\leqq i-\sum_{i=0}^{n-1} n_i\leqq\sum_{i=0}^{k+1}(n_i-1)\,\, \textrm{and} \\

& \hspace{0.7cm} \sum_{i=0}^{l}(n_i-1)+1\leqq j-\sum_{i=0}^{n-1} n_i\leqq\sum_{i=0}^{l+1}(n_i-1),  \\

\end{array}
\right.
$$
\noindent for some $k,\,l\in\{1,2,\ldots, n-2\}$. Here, $n_0=0$, and indexes $I_i$ and $J_i$, $i=1,2,3,4$, are the following:
$$
\left\{
\begin{array}{@{\,}ll}
\medskip

I_1 =k+1, i-\sum_{i=1}^k n_i \, ;\, J_1 =l+1, i-\sum_{i=1}^l n_i \, ; \\ 
\medskip

I_2 =k+1, i-\sum_{i=1}^k n_i \, ;\, J_2 =l+1, j-\sum_{i=1}^l n_i -\sum_{i=1}^l (n_i-1) \, ; \\ 
\medskip

I_3 =k+1, j-\sum_{i=1}^k n_i -\sum_{i=1}^k (n_i-1) \, ;\, J_2 =l+1, i-\sum_{i=1}^l n_i \, ; \\ 
\medskip

I_4 =k+1, j-\sum_{i=1}^k n_i -\sum_{i=1}^k (n_i-1) \, ;\, J_2 =l+1, j-\sum_{i=1}^l n_i -\sum_{i=1}^l (n_i-1). \\ 
\end{array}
\right.
$$
\end{proof}

We call $h_{i,j}^l-h_{j,i}^l$ the alternation of $h_{i,j}^l$, denoted by 
$h_{i,j}^{l,A}$, where $i<j$, similarly, We call $p_{i,j}^{l,m}-p_{j,i}^{m,l}$ the alternation of $p_{i,j}^{l,m}$, denoted by $p_{i,j}^{l,m,A}$, where $i<j$.

\begin{theorem} \label{ncbp} Let $\{\, b_{l,1}, b_{l,2},\ldots, b_{l,n_l-1}\,\, |\,\, l=1,2,\ldots, n-1 \,\}$ be a set of oriented loops as in the proof of Theorem \ref{ncbm} and let $H_l=(h_{i,j}^l)$, $P_{l,m}=(p_{i,j}^{l,m})$, $1\leqq l\neq m\leqq n-1$. Then the following properties hold: 
\begin{enumerate}
\item[$(I)$] $|h_{i,j}^{A, l}|=0$ or $1$. 
\item[$(II)$] $|p^{l,m,A}_{i,j}|=0$ or $1$.
\end{enumerate}
\end{theorem}

\begin{proof} Note that $h_{i,j}^l=lk(b_{l,i}, b_{l,j}^{+})$ and $p_{i,j}^{l,m}=lk(b_{l,i}, b_{m,j}^{+})$, and that if $D_l\cap b_{l,i}\cap b_{l,i+1}\neq\emptyset$, then $i(b_{l,i}\cap D_l, b_{l,i+1}\cap D_l)=-1$ and if $D_n\cap b_{m,j}\cap b_{m,j+1}\neq\emptyset$, then $i(b_{m,j}\cap D_n, b_{m,j+1}\cap D_n)=1$, where $i\in\{1, 2, \ldots , n_l-2\}$, $j\in\{1, 2, \ldots , n_m-2\}$, and $l, m\in\{1, 2, \ldots , n-1\}$.
\medskip

First, we show $(I)$ in the case $j=i+1$, $i\in\{1, 2, \ldots , n_l-2\}$, $l\in\{1, 2, \ldots , n-1\}$.

If $b_{l,i}\cap b_{l,i+1}=\emptyset$, i.e., $b_{l,i}\cap b_{l,i+1}\cap D_l=b_{l,i}\cap b_{l,i+1}\cap D_n=\emptyset$, then $h_{i,i+1}^{l,A}=lk(b_{l,i}, b_{l,i+1}^+)-lk(b_{l,i+1}, b_{l,i}^+)=0$. 

If $b_{l,i}\cap b_{l,i+1}\cap D_l=\emptyset$ and $b_{l,i}\cap b_{l,i+1}\cap D_n$ consists of one point, then $h_{i,i+1}^{l,A}=1$.

If $b_{l,i}\cap b_{l,i+1}\cap D_l$ consists of one point and $b_{l,i}\cap b_{l,i+1}\cap D_n=\emptyset$, then $h_{i,i+1}^{l,A}=-1$.

If both sets $b_{l,i}\cap b_{l,i+1}\cap D_l$ and $b_{l,i}\cap b_{l,i+1}\cap D_n$ consist of one point, then $h_{i,i+1}^{l,A}=0$. 
\smallskip

Next, we show $(I)$ in the case $j\neq i+1$, $i\in\{1, 2, \ldots , n_l-1\}$, $l\in\{1, 2, \ldots , n-1\}$.

If $b_{l,i}\cap b_{l,i+1}=\emptyset$, i.e., $b_{l,i}\cap b_{l,i+1}\cap D_l=b_{l,i}\cap b_{l,i+1}\cap D_n=\emptyset$, then $h_{i,i+1}^{l,A}=0$. 

If $b_{l,i}\cap b_{l,i+1}\cap D_l=\emptyset$ and $b_{l,i}\cap b_{l,i+1}\cap D_n$ consists of one point, then $h_{i,i+1}^{l,A}=1$ or $-1$.

If $b_{l,i}\cap b_{l,i+1}\cap D_l$ consists of one point and $b_{l,i}\cap b_{l,i+1}\cap D_n=\emptyset$, then $h_{i,i+1}^{l,A}=1$ or $-1$.

If both sets $b_{l,i}\cap b_{l,i+1}\cap D_l$ and $b_{l,i}\cap b_{l,i+1}\cap D_n$ consist of one point, then $h_{i,i+1}^{l,A}=0$. Therefore, $(I)$ holds. 
\medskip

Note that $b_{l,i}\cap b_{m,j}\cap D_l=\emptyset$ and $b_{l,i}\cap b_{m,j}\cap D_n$ consists of at most one point for any $l\neq m$, $l,m\in\{1, 2, \ldots , n-1\}$, $i\in\{1, 2, \ldots , n_l-1\}$, $j\in\{1, 2, \ldots , n_m-1\}$. 

Next, we show $(II)$. Let $l\neq m$, $l$, $m\in\{1, 2, \ldots , n-1\}$, $i\in\{1, 2, \ldots , n_l-1\}$, $j\in\{1, 2, \ldots , n_m-1\}$.

If $b_{l,i}\cap b_{m,j}\cap D_n=\emptyset$, then $p_{i,j}^{l,m,A}=0$. 

If $b_{l,i}\cap b_{m,j}\cap D_n$ consists of one point, then $p_{i,j}^{l,m,A}=1$ or $-1$. Therefore $(II)$ holds. 
\smallskip

Thus, $(I)$ and $(II)$ holds. The proof is completed. 
\end{proof}

\begin{corollary} Let $L$ be an $n$-component Brunnian link, and let $C_L=\cup _{i=1}^n D^2_i$ be a C-complex for $L$ obtained by \lq\lq Algorithm\rq\rq\, in section 3. If clasp singularities in $D_i$ arrange as in Figure  for each $i\in\{1, 2, \ldots , n-1\}$, then there is a Seifert matrix $M_L$ for $L$ of the form as in Theorem 5 so that the following equalities hold for some $s_2^l\in{\mathcal S}(2)$ and $s_j^l\in{\mathcal S}(s_2^l, s_3^l, \ldots , s_{j-1}^l, j)$, $j=3, 4, \ldots , n_l-1$: 
$${\mathcal G}_2\bigl([s_2^l]\bigr)=\bigl(h_{1,2}^{l,A}\bigr), 
$$ 
$${\mathcal G}_{s_2^l, s_3^l, \ldots , s_{j-2}^l, j-1}\bigl([s_{j-1}^l]\bigr)=\bigl( h_{1,j-1}^{l,A}, h_{2,j-1}^{l,A}, \ldots , h_{j-2,j-1}^{l,A}\bigr), 
$$ 
and 
$$\left|p^{l,m, A}_{i,j}\right|=0 \,\,{\textrm or} \,\,\,\, 1.
$$
\end{corollary}

\begin{proof} This is an immediate consequence of Theorem \ref{ncbp}.
\end{proof}

\begin{corollary} \label{co3} Let $M_L$ be a Seifert matrix for an $n$-component Brunnian link $L$, $n\geqq 3$ as in Theorem 5. Then the following equality holds: $M_L-tM_L^T$
$$
=\left(
\begin{array}{c|c|c|c|c|c|c|c}
E_{n_1}^* & 0 & 0 & 0 & F_{n_1} & 0 & 0 & 0 \\ \hline
0 & E_{n_2}^* & 0 & 0 & 0 & F_{n_2} & 0 & 0 \\ \hline
0 & 0 & \ddots & 0 & 0 & 0 & \ddots  & 0 \\ \hline
0 & 0 & 0 &  E_{n_{n-1}}^* & 0 & 0 & 0 &  F_{n_{n-1}} \\ \hline
-tF_{n_1}^T & 0 & 0 & 0 & H_1^* & P_{1,2}^* & \cdots & P_{1,{n-1}}^* \\ \hline
0 & -tF_{n_2}^T & 0 & 0 & P_{2,1}^* & \ddots & \ddots & \vdots \\ \hline 
0 & 0 & \ddots & 0 & \vdots & \ddots & H_{n-2}^* & P_{{n-2},{n-1}}^* \\ \hline 
0 & 0 & 0 & -tF_{n_{n-1}}^T & P_{{n-1},1}^* & \cdots & P_{{n-1},{n-2}}^* & H_{n-1}^* \\
\end{array}
\right).
$$
\noindent Here, $E_{n_l}^*\equiv (1-t)E_{n_l}$, $H_l^*\equiv H_l-tH_l^T$, and $P_{l,m}^*\equiv P_{l,m}-tP_{l,m}^T$. And then, all of the determinants of the following proper submatrices of $M_L-tM_L^T$ equal zero:
$$
\left(
\begin{array}{c|c|c|c|c|c|c|c}
E_{n_{j_1}}^* & 0 & 0 & 0 & F_{n_{j_1}} & 0 & 0 & 0 \\ \hline
0 & E_{n_{j_2}}^* & 0 & 0 & 0 & F_{n_{j_2}} & 0 & 0 \\ \hline
0 & 0 & \ddots & 0 & 0 & 0 & \ddots  & 0 \\ \hline
0 & 0 & 0 &  E_{n_{j_k}}^* & 0 & 0 & 0 &  F_{n_{j_k}} \\ \hline
-tF_{n_{j_1}}^T & 0 & 0 & 0 & H_{j_1}^* & P_{j_1,j_2}^* & \cdots & P_{j_1,j_k}^* \\ \hline
0 & -tF_{n_{j_2}}^T & 0 & 0 & P_{j_2,j_1}^* & \ddots & \ddots & \vdots \\ \hline 
0 & 0 & \ddots & 0 & \vdots & \ddots & H_{j_{k-1}}^* & P_{j_{k-1},j_k}^* \\ \hline 
0 & 0 & 0 & -tF_{n_{j_k}}^T & P_{{j_k},{j_1}}^* & \cdots & P_{j_k,j_{k-1}}^* & H_{j_k}^* \\
\end{array}
\right),
$$

\noindent where $k\in\{1, 2, \ldots, n-2\}$ and $j_1,j_2\ldots, j_{k}\in\{1, 2, \ldots, n-1\}$.
\end{corollary}

\begin{proof} The former statement is an immediate consequence of Theorem \ref{ncbm}. We shall show the latter statement. Let $L=\cup_{i=1}^n K_i$. By the definition of Brunnian properties, any proper sublink $\cup_{j\in N'} K_j$ of $L$ is unlink, where $N'$ is a proper subset of $\{1, 2, \ldots, n\}$. Therefore, we see that for any $N'$, the equality $\Delta_{\cup_{j\in N'}K_j}(t)=0$ holds, i.e., the latter statement of this theorem holds. 
\end{proof}

%\begin{remark} Let $H_l^*=(h_{i,j}^{l*})$, $P_{l,m}^*=(p_{i,j}^{l,m*})$. Then t%he following conditions hold: if $i=j$, then $h_{i,j}^{l*}=h_{i,j}^l(1-t)$ and %$p_{i,j}^{l,m*}=p_{i,j}^{l,m}(1-t)$; otherwise, then $h_{i,j}^{l*}=h_{i,j}^l(1-%t)$ or $h_{i,j}^l(1-t)\pm t$, and $p_{i,j}^{l,m*}=p_{i,j}^{l,m}(1-t)$ or $p_{i,%j}^{l,m}(1-t)+t=p_{j,i}^{m,l}(1-t)+1$. 
%\end{remark}

\bigskip

\section{acknowledgements}
\medskip

The author would like to appreciate Professor Kunio Murasugi and Professor Shinji Fukuhara for some information and discussions. She would also thank Professor Teruhisa Kadokami for his useful advices and thank Professor Kenichi Tamano for his kind encouragements.

\end{document}